\documentclass[11pt]{article}
\usepackage{latexsym}
\usepackage{amsfonts}
\usepackage{amsmath}
\makeatletter

\def\@footnotetext#1{\insert\footins{%

\footnotesize 
 
 \interlinepenalty\interfootnotelinepenalty

 \splittopskip\footnotesep

 \splitmaxdepth \dp\strutbox \floatingpenalty \@MM

 \hsize\columnwidth \@parboxrestore

 \edef\@currentlabel{\csname 
p@footnote\endcsname\@thefnmark}\@makefntext
 {\rule{\z@}{\footnotesep}\ignorespaces
#1\strut}}}

\def\abstract{\small\quotation{\hskip-\parindent\sc Abstract.}}
\def\classification{\@ifnextchar [{\@xfootnotenext}%
 {\begingroup\let\protect\noexpand
 \xdef\@thefnmark{}\endgroup
 \@footnotetext}}

\title {}

\begin{document}

\classification {{\it 2000 Mathematics
Subject Classification:} Primary 68Q17; Secondary 20F10, 20P05.}

\begin{center}

{\bf \Large  Assessing security of some 
   group based cryptosystems}

\bigskip

{\bf Vladimir Shpilrain }

\end{center}

\medskip

\begin{abstract}
\noindent    One of the possible generalizations of the 
{\it discrete logarithm problem} to arbitrary groups is the so-called 
 {\it  conjugacy search  problem} (sometimes erroneously called just 
the {\it conjugacy problem}) : given two elements $a, b$ of a group $G$ 
 and the information that $a^x=b$ 
for some $x \in G$, find at least one particular element $x$ like that. Here 
$a^x$ stands for $xax^{-1}$. 
The computational difficulty of this  problem 
in some particular groups has been  used in several group based 
cryptosystems. Recently, a few preprints have been 
in circulation that suggested various ``neighbourhood search" type
heuristic attacks on the conjugacy search   
problem. The goal of the present survey is to stress a (probably well known) 
fact that these heuristic attacks alone are not a threat to the security of 
a cryptosystem, and, more importantly, to suggest a more credible 
approach to assessing security of group based cryptosystems. Such an 
approach should be necessarily based on the concept of the {\it average 
case complexity} (or {\it expected running time}) of an algorithm.

 These arguments support the following conclusion: although it is 
generally feasible to base the security of a cryptosystem on the 
difficulty of the conjugacy search   problem, the group $G$ itself 
(the ``platform") has to be chosen very carefully. In particular, 
experimental as well as theoretical evidence collected so far 
makes it appear likely that braid groups are {\it not} a good choice 
for the platform. We also reflect on possible replacements.
\end{abstract} 

\bigskip

\section{Introduction}

 Let $G$ be a group. The {\it conjugacy problem} (or {\it conjugacy decision 
problem}) for $G$ is: 
given two elements $a, b \in G$, find out whether  
there is $x \in G$ such that $a^x=b$, where $a^x$ stands for $xax^{-1}$.

 On the other hand, the  {\it conjugacy search problem} 
(or {\it witness conjugacy  problem}) is:
given two elements $a, b \in G$   and the information that $a^x=b$ 
for some $x \in G$, find at least one particular element $x$ like that.

 The conjugacy problem is of interest in group theory, but of no 
interest in complexity theory. In contrast, the  conjugacy search   problem
is of interest in complexity theory, but of no 
interest in group theory. Indeed, if you know that $a$ is conjugate to 
$b$, you can just go over elements of the form $a^x$ and compare them to 
$b$ one at a time, until you get a match. (We implicitly use here an 
obvious fact that a group with solvable conjugacy problem also has 
solvable word problem.) This straightforward  algorithm is at least 
exponential-time in the length of $b$, and  therefore is considered   
infeasible for any practical purposes. 

 Thus, if no other algorithm is known for the  conjugacy search   problem 
in a group $G$, it is not unreasonable to claim that $x \to a^x$ is a 
one-way function and build a (public-key) cryptosystem on that. 
A simple way of doing it was suggested in \cite{KLCHKP}, and a more 
sophisticated and more general one in the seminal paper \cite{AAG} (see 
also \cite{AAFG}). 
 Just to make our exposition self-contained, we give a brief exposition 
of both protocols in Section \ref{protocols}. 
In both cryptosystems a braid group 
$B_n$ is used as the ``platform". We are not going into any technical 
details about braid groups here since they are irrelevant to the subject of 
the present survey. It is sufficient to say that the choice of a braid group 
as the platform is highly questionable from the security point of view,  
although quite understandable from the marketing point of view. We discuss
 this in more detail in Section \ref{vs}. 

There are numerous deterministic algorithms for solving the conjugacy search problem
 in braid groups 
(see e.g. \cite{BKL1},  \cite{BKL2}, \cite{EM}, \cite{GM}) 
that  are believed to be polynomial-time. As soon as one of them  
is {\it  proved} to be polynomial-time,  braid groups 
are  out of the game. We are going to explain in Section \ref{strike} 
 why the existence of a polynomial-time deterministic algorithm 
``usually" implies the existence of a deterministic algorithm that has 
linear-time complexity {\it on average}.

 Several ``neighbourhood search" type (in a group-theoretic context sometimes 
also  called  ``length based") and other  heuristic algorithms proposed recently 
\cite{GKTTV}, \cite{Gebhardt}, \cite{HS}, \cite{HT}, \cite{LeePark} \cite{Ushakov} 
show very fast performance on ``most" inputs which 
has driven some of the authors to a somewhat premature conclusion 
that cryptosystems \cite{AAG} and \cite{KLCHKP} are insecure if a braid group 
is chosen as the platform. Whereas we agree with the conclusion itself, we have 
to be fair and admit that at the time of this writing, 
there is not enough evidence  to rigorously prove it. 

 The reason is as follows: fast performance of  heuristic algorithms 
on ``most" inputs means that {\it generic-case} time-complexity of those 
algorithms is low (probably linear-time). However there is a ``small" 
set of ``noncooperative" inputs on which a given heuristic 
algorithm may work very slowly or not terminate at all. Therefore, a 
countermeasure against 
 an heuristic  attack like that would be just using several rounds of the 
same protocol 
with different random choices of elements. (This will affect the 
efficiency of the protocol, but we are assuming that the original 
protocol is efficient enough so that repeating it a few dozen times 
does not make it impractical.) 
If, for example, the probability
of success of a particular heuristic algorithm in a single round is 
$0.9$, then the probability of success in, say, 50 rounds is $0.9^{50}$, 
which is as good as 0. Using several runs of the same ``neighbourhood search" type 
heuristic algorithm on the same input will not help, in general, to improve the 
probability of success since being ``noncooperative" (with respect to a given 
algorithm) is usually an intrinsic property of an input. 

\medskip 

\noindent  {\bf Remark.} It might be tempting to say that   
 narrowing down  the set of inputs to a subset of ``noncooperative" 
inputs for a given heuristic algorithm will render this particular 
algorithm ineffective (or, at least, less effective). 
Although this countermeasure is theoretically feasible, 
it will not work in real life because in real life, 
once a cryptosystem has entered  service, its parameters cannot be 
significantly changed. In particular, the pool from which inputs are 
(randomly) selected cannot be adjusted in response to 
one attack  or another. 
This means a potential   attacker has what looks like an 
unfair advantage, but these are rules of the game.

\medskip 

 To properly take into account the set of ``noncooperative" inputs, one has 
to use the notion of the {\it average-case} time-complexity 
(or {\it expected running time}) of an algorithm rather 
than  generic-case. Differences between these properties for various 
algorithms in group theory have been studied in-depth in \cite{KMSS1}
 and \cite{KMSS2}.  We give more details in Section 2; here we just mention 
that for the average-case complexity of an algorithm to be defined, 
this algorithm should terminate for any input. 

 If one wants to enhance the average-case performance of a  
deterministic algorithm ${\cal D}$, 
one can run a (relatively fast) heuristic algorithm ${\cal H}$ in parallel 
with   ${\cal D}$. This composite algorithm ${\cal H} \mid \mid {\cal D}$ 
is still deterministic, and its 
 average-case complexity is therefore defined. 
If the average-case complexity of ${\cal H} \mid \mid {\cal D}$  can be 
computed (or estimated) 
and turns out to be low (say, linear- or quadratic-time), thenn\footnote{This is 
the recommended notation for ``then and only then".}
 the relevant cryptosystem can be considered in serious danger.  

 All this is probably well-known (at least, on the intuitive level). 
However, there is a spectacular result in \cite{KMSS2} that may come as 
a surprise. Adopted in our situation, it says that, if the algorithm 
${\cal H}$ has linear-time strong  generic-case complexity, and ${\cal D}$ has a  
subexponential complexity, then the composite 
algorithm ${\cal H} \mid \mid {\cal D}$ has linear-time average-case complexity. 
Therefore, a cryptosystem whose security is 
based on the conjugacy search   problem, can be considered reasonably 
secure only if there is no deterministic algorithm that solves this problem 
in the platform group in subexponential time. 
 We discuss this in more detail 
 in Sections \ref{complexity} and  \ref{strike}. 

 Finally, in the concluding Section \ref{group?}, we discuss possible 
replacements for braid groups in the otherwise promising cryptosystem 
\cite{AAG}.

\section{Cryptographic protocols involving the conjugacy search  problem}
\label{protocols}

 Let $G$ be a group with  solvable word problem. Recall that for 
$a, x \in G$, the notation $a^x$ stands for $xax^{-1}$. 

 We start with a simpler protocol, due to Ko et. al. \cite{KLCHKP}. 
\medskip 

\noindent  {\bf (0)}  An  element $a \in G$ is published. 

\noindent  {\bf (1)} Alice picks a private $x \in G$ and sends $a^x$ to Bob. 

\noindent  {\bf (2)} Bob picks a private $y \in G$ and sends $a^y$ to Alice. 

\noindent  {\bf (3)} Alice computes $(a^y)^x = a^{yx}$, and Bob computes $(a^x)^y = a^{xy}$.
\medskip 

 If $x$ and $y$ are chosen from a pool of commuting elements of the group 
$G$, then $xy=yx$, and therefore, Alice and Bob get a common private key $a^{xy}=a^{yx}$. 
In case $G=B_n$, a braid group, the pool of commuting elements is quite large, 
so there is no problem with setting up a protocol like that if a braid group 
is selected as the platform. 
\medskip 

 Another protocol, due to Anshel et. al. \cite{AAG}, is more complex, 
but it is more general in the sense that there are no requirements  
on the group $G$ other 
than to have  solvable word problem. This really makes a difference 
and gives a great advantage to the protocol of \cite{AAG} over that 
of \cite{KLCHKP}. When the braid group platform is rigorously 
shown to be insecure (conceivably, this will happen rather soon), 
it will be much 
easier to repair the protocol of \cite{AAG} (by choosing a different 
platform) than that of \cite{KLCHKP}. In fact,  an attempt  
at such repairing has already been made in \cite{AAFG}, but it 
definitely should be taken further (away from braid groups). 
 
\medskip 

\noindent  {\bf (0)} Elements $a_1,...,a_k, b_1,...,b_m \in G$ are published. 

\noindent  {\bf (1)} Alice picks a private $x \in G$ as a word in 
$b_1,...,b_m$ (i.e.,  $x=x(b_1,...,b_m)$)  and sends 
$a_1^x,...,a_k^x$ to Bob. 

\noindent  {\bf (2)} Bob picks a private $y \in G$ 
as a word in $a_1,...,a_k$ and sends 
$b_1^y,...,b_m^y$ to Alice. 

\noindent  {\bf (3)} Alice computes $x(b_1^y,...,b_m^y) = x^y = yxy^{-1}$, 
and Bob computes $y(a_1^x,...,a_k^x) = y^x = xyx^{-1}$. 

\noindent  {\bf (4)} Alice and Bob come up with a common private key $xyx^{-1}y^{-1}$ 
(called the \emph{commutator} of $x$ and $y$) as follows: Bob just multiplies 
$xyx^{-1}$ by $y^{-1}$ on the right, while Alice multiplies $yxy^{-1}$ by 
$x^{-1}$ on the right, and then takes the inverse of the whole thing: 
$(yxy^{-1}x^{-1})^{-1} = xyx^{-1}y^{-1}$. 
\medskip 

  Finally, we re-iterate the  point that the security of both 
protocols described in this section apparently 
rests on the computational difficulty 
of the  conjugacy search   problem in the platform group $G$.

\section{Generic-case vs. average-case complexity} 
\label{complexity} 

 In this section, we offer an informal discussion of generic- and 
average-case complexity of algorithms, referring to \cite{KMSS1}
 and \cite{KMSS2} for formal definitions and results. 

To   discuss   \emph{generic-case} complexity, which deals with the
performance of an algorithm $\Omega$ on ``most'' inputs, we  first need
a notion of which sets of inputs  are \emph{generic}. Let $\nu$ be 
 an arbitrary additive function with values in 
$[0,1]$ defined on some 
subsets of the set $X$ of possible inputs.
A subset $T \subset X$ is called \emph{generic with respect
to}  $\nu$ (or just $\nu$-generic) 
if $\nu(X - T) = 0$.  Then, for example,  we would
say that   $\Omega$ has \emph{polynomial-time generic-case  complexity 
with respect to $\nu$} if $\Omega$ runs in polynomial time on all inputs 
from some  subset $T$ of  $X$  which is generic
with respect to $\nu$.
Of course, one can  define   generic-case complexity being in any
 complexity class $\mathcal{C}$, not just polynomial-time.

Generic-case complexity therefore does not take into account behaviour 
of an algorithm on the set $X - T$ of ``noncooperative" inputs. In fact, 
on some of those inputs, the algorithm in question may not terminate 
at all. In contrast, the average-case complexity of an algorithm 
takes into account all possible inputs; in particular, for the average-case 
complexity of an algorithm to be defined, this algorithm should terminate 
for any input.

Here we have to assume that there is some kind of complexity $|w|$ 
defined for all elements $w$ of the set $X$ of possible inputs. We will  call it 
the ``length" of $w$ to simplify the language. 

We will say that a (discrete) probability measure $\mu$ is \emph{length-invariant} if for any
elements  $w, w'\in X$ with $|w|=|w'|$ one has $\mu(w)=\mu(w')$.
Requiring that a measure be length-invariant is a very natural
assumption since most algorithm complexity classes are defined in terms of
the length of an input.

 Now let $\cal P$ be a property that elements of $X$ may or may not have, 
and let $\cal A$ be an algorithm
which for every $w \in X$ decides, in time $T(w) < \infty$, whether or 
not $w$ has the property $\cal P$. 

Let $f(n)$ be a non-decreasing positive function. We say that 
 $\cal A$ has
average case time-complexity bounded by $f(n)$ relative to $\mu$ if

\[
\int_{X} \frac{T(w)}{f(|w|)} \mu(w) =\sum_{w\in X}
\frac{T(w)}{f(|w|)} \mu(w)  < \infty.
\]

 Note that if this is the case, then the ``$\mu$-expected running time"
$\sum_{|w|=n}  T(w) \mu(w)$  for inputs of length 
$n$ is $o(f(n))$.

 This definition may look somewhat scary, but this is just because of 
the probability measure $\mu$. We note that to define a meaningful probability measure 
on an infinite group is a very non-trivial theoretical problem (see \cite{BMS}).
 However, for practical purposes, e.g. for applications in cryptography, 
the theoretical mumbo-jumbo can be avoided because an infinite group in this context 
 can be considered as just a big finite group, and the following kind of 
practical experiment can be conducted. 

\medskip 

\noindent  {\bf Experiment.} Usually, ``neighbourhood search" type 
heuristic algorithms show very fast performance on 
``most" inputs, i.e., they tend to have very low generic-case complexity. 
In particular, this is usually the case with {\it genetic algorithms} 
(see e.g. \cite{Ushakov}); 
 these are subtle combinations of ``neighbourhood search" and Monte Carlo type 
algorithms. 

 Experiments conducted in \cite{GKTTV, HS, Hughes, LeePark} show a rather 
high rate of success 
(80\% and up) of some particular heuristic algorithms. However, as we 
have explained, this alone is not enough to make any conclusions 
about the security of a relevant cryptosystem. The ``right" experiment 
should be arranged as follows. 

 Run your heuristic algorithm ${\cal H}$ in parallel with a  
deterministic algorithm ${\cal D}$ of your choice. Call this combined 
 algorithm ${\cal A}$. Now suppose the algorithm ${\cal H}$  terminates   
in at most $h(n)$ seconds on $b$\% of randomly selected inputs 
of length  $n$.  
On the remaining (100-$b$)\% of those inputs, the 
deterministic algorithm ${\cal D}$  kicks in and terminates in 
at most $d(n)$ seconds. Then the expected running time of the 
algorithm ${\cal A}$ on inputs of length  $n$ is at most  
$e(n) = \frac{1}{100} \cdot [h(n) \cdot b + d(n) \cdot (100-b)]$ seconds. 
 Then you try to extrapolate the  function $e(n)$  by a polynomial 
of a small degree with ``reasonably small" coefficients. If you succeed, 
then you can  claim that the relevant cryptosystem is likely  to
 be insecure. $\Box$ 

\medskip 

\noindent  {\bf Remark.} Once again, we emphasize the point that 
 generating only those inputs that would be non-cooperative for 
${\cal H}$ is possible in theory but not in  
   real life  because in  real  life, 
once a cryptosystem has entered  service, its 
 parameters cannot be significantly changed. Therefore, an   attacker 
clearly has an  edge here because he can choose  ${\cal H}$ based 
on the parameters of a cryptosystem, and not the other way around. 

\medskip

For the next section, we 
 need the notion of \emph{strong generic-case} complexity.
 For a subset $S \subseteq X$, let $S_n = \{w \in S, ~|w|=n\}$. 
Let $\nu$ be  an additive  function with values in $[0,1]$ defined on some 
subsets of the set $X$ of possible inputs. We assume that 
if $\nu(S)$ is defined, then $\nu(S_n)$ is defined for any $n$. 
We call $S$ \emph{strongly $\nu$-generic (or exponentially $\nu$-generic)}   
if $\nu(X_n - S_n) \to 0$ exponentially fast 
as $n \to \infty$. Accordingly, an algorithm $\Omega$ has, say,  
\emph{polynomial-time strong  generic-case  complexity 
with respect to $\nu$} if $\Omega$ runs in polynomial time on all inputs 
from some  subset $S$ of  $X$  which is strongly $\nu$-generic.

\section{``If the strike is not settled quickly, it may last a while"}
\label{strike}

 This was one of the top 10 ``Funniest headlines of the year" in 1995. 
However, it appears that this particular headline may not be that 
funny after all, but rather can be an instance of a universal rule, something 
similar (in magnitude) to the ``coupling principle" in quantum physics 
in its effect on our lives. Probably everyone has been in a situation 
where he/she receives a message or a letter that requires a response; then, 
if you do not respond immediately, you forget about the whole thing 
and respond much later, usually after a reminder. Or, you start a project 
(say, a paper), and you either complete it in a few weeks or it gets 
stretched to several months (or years). There are numerous 
other real-life  manifestations of the same philosophical principle, 
including the one alluded to in the aforementioned newspaper headline. 
The following result from \cite{KMSS2} (Proposition 3.2) might give a 
theoretical support 
to this intuitive principle. We give a simplified version here that better 
fits in with our subject. 

\medskip

\noindent  {\bf Theorem.} Let $\mu$ be an arbitrary length-invariant discrete probability 
measure. Suppose there is an algorithm ${\cal H}$, deterministic or 
not, that solves a given problem ${\cal P}$ strongly $\mu$-generically in linear 
(or quadratic) time with respect to the complexity of an input. Suppose also 
that there is a
subexponential-time deterministic algorithm ${\cal D}$ for solving  ${\cal P}$. 
Then the   (deterministic) algorithm ${\cal H} \mid \mid {\cal D}$  solves  
${\cal P}$ in linear (resp. quadratic) time on average relative to  $\mu$.
 
\medskip

 Informally speaking, this result says that, for a reasonably natural problem,
 either there  is a {\it very} fast (on average) algorithm 
for solving it,   or there is no fast (on average) algorithm at all. 
Or, yet in other words, there is a gap between ``very fast" 
and   ``slow". 

  We note in passing (although it is not directly related to the 
subject of this survey) that the conjugacy  problem for braid groups 
is solvable  generically in linear time by \cite[Theorem C]{KMSS1}. 
We do not know however whether it is solvable {\it strongly} 
generically in linear time. Neither do we have any useful results on 
the generic-case complexity of the {\it conjugacy search   problem} 
for braid groups.

\section{Marketability vs. security} 
\label{vs}

 It is a fact that abstract groups, unlike numbers, are not 
something that most people learn at  school. There is therefore an obvious 
communication problem involved in marketing a cryptographic  product 
that uses abstract groups one way or another. Braid groups clearly 
have an edge here because to explain what they are, one can draw 
simple pictures, thus alleviating the fear of the unknown.   
The fact that braid groups cut across  many different areas of mathematics 
(and physics) helps, too, since this gives more credibility  
to the hardness of the relevant problem (the  conjugacy search   problem 
in our case). 

 We can recall that, for example, confidence in the 
 security of the RSA cryptosystem is based on literally centuries-long 
history of attempts by thousands of people, including such authorities 
as Euler and Gauss, at factoring integers fast. The history of braid groups 
goes back to 1927, and again, thousands (well, maybe hundreds) of people, 
including prominent mathematicians like Artin, Thurston, V. F. R. Jones, 
and others have been working on various aspects, including algorithmic ones,
 of these groups. 

 On the other hand, from the security point of view, the fact that 
braid groups cut across so many different areas can be a disadvantage, 
because different areas provide different  tools for solving  a problem at
 hand (in our situation, the conjugacy search   problem). Indeed, people 
from several different areas (group theory, topology, combinatorics) 
have already contributed to the solution of this problem. This increases  
the  odds of finding a subexponential-time deterministic algorithm, 
and therefore breaking the relevant cryptosystem (see Section \ref{strike}). 
 Furthermore, braid groups turned out to be linear \cite{Bigelow}, 
\cite{Krammer}, which makes them potentially vulnerable to linear 
algebraic attacks (see e.g.  \cite{Hughes}, \cite{LeeLee}), and this alone 
is a serious security hazard.

 The pioneering paper \cite{AAG} has brought combinatorial group theory 
into cryptography in a very serious and promising way. The choice of 
braid groups as the platform was probably inevitable at that time, 
for the reasons outlined above. At the same time, this  has  paved the way for 
 engaging other groups more easily in the future. It is probably time 
now to take advantage  of this opportunity and start a serious search for a 
secure platform rather than keep chewing on braid groups. This search is not 
going to be easy, as we try to explain in the next section.

\section{Search for a platform} 
\label{group?}

 Before we can start a search, we have to put down the properties 
that we want from a group $G$ in this context.

 It seems reasonable to start with a property 
which, although not mathematical, appears to be mandatory 
if we want our cryptographic product to  be used in real life (see the 
previous section): 
\smallskip 

\noindent {\bf (P0)} The group has to be well known. More precisely, 
the conjugacy search   problem in the group either 
has to be well studied or can be reduced to a well known problem 
(perhaps, in some other area of mathematics). 
\smallskip 

 We note in passing that this property
already narrows down the list of candidates quite a bit. 

 The following two are mathematical properties. 

\smallskip 

\noindent {\bf (P1)} The word problem in $G$ should have 
a fast (linear- or quadratic-time) solution by a deterministic algorithm.
\smallskip 

 This is required for an efficient common key extraction by legitimate 
parties.
\smallskip 

\noindent {\bf (P2)} The conjugacy search   problem should {\it not} have 
a subexponential-time solution by a deterministic algorithm. 
\smallskip 

 We point out here that {\it proving} a group to have  (P2) should be 
extremely difficult, if not impossible. This is, literally, a million-dollar 
problem (see \cite{Clay}). The  property (P2) should be therefore considered in 
conjunction with (P0), i.e., the only realistic 
evidence of a group $G$ having the property 
(P2) can be the fact that sufficiently many people have been studying 
the  conjugacy search   problem  in $G$ over sufficiently long time. 
\smallskip 

 The last property is somewhat informal, but it is rather important 
for practical implementations: 
\smallskip 

\noindent {\bf (P3)} There should be a way to disguise elements of $G$ 
so that it would be  impossible to recover  $x$ from $xax^{-1}$ just by
inspection.  In particular, if $G$ is given by means of generators 
and relations, then at least some of these relations should be very  
short. 
\smallskip 

 We see now that braid groups have (P0), (P1), (P3), but (most likely) 
not (P2). Not having (P2) is, of course, a grave security 
hazard, so as soon as it is  proved  that braid groups do not have (P2), 
they are  out of the game. 

 There are groups that have  (P1), most likely have (P2), and to a 
reasonable extent have (P3).  These are groups with solvable word problem, 
but unsolvable conjugacy problem (see e.g. \cite{Miller}).
However, groups like that tend not to have  the property (P0) because, 
as we have mentioned in the Introduction, the conjugacy search   problem is of 
no independent interest in group theory. Group theorists therefore did not bother 
to study the conjugacy search   problem in these groups once it had been 
proved that the conjugacy problem is algorithmically unsolvable. 
It would probably make sense now to reconsider these groups.

 Another possible strategy would be to base the search on braid groups, i.e., 
to ``distort" braid groups one way or another (e.g. to discard a couple 
of defining relations)  to get groups with the property 
(P2), while keeping (P1) intact. 
The  property  (P0) will be lost, of course, but the new 
groups   might have marketing potential nonetheless: ``Look, these 
are like braid groups, only better". 

 We do not have any more specific suggestions at this time.

\noindent 
 Department of Mathematics, The City  College  of New York, New York, 
NY 10031 
\smallskip

\noindent {\it e-mail address\/}: 
shpil@groups.sci.ccny.cuny.edu  

\smallskip

\noindent {\it http://www.sci.ccny.cuny.edu/\~\/shpil/} \\

\end{document}